\documentclass[11pt]{amsart}
\usepackage{epsfig}
\usepackage{verbatim}
\usepackage{amssymb}
\usepackage{amsfonts}
\usepackage{txfonts}
\usepackage{graphicx}
\newcommand{\avg}{\mathrm{avg}}
\newcommand{\var}{\mathrm{var}}
\newcommand{\clus}{\mathrm{clus}}

\newcommand{\bvec}{\left[ \begin{array}{c}}
\newcommand{\evec}{\end{array} \right]}
\newcommand{\refl}{\mathrm{refl}}

\newcommand{\eproof}{\hfill\rule{2.2mm}{3.0mm}}
\newcommand{\esubproof}{\hfill$\Box$}
\newcommand{\Proof}{\noindent {\bf Proof.~~}}

\newcommand{\Z}{{\mathbb Z}}

\newcommand{\Q}{{\mathbb Q}}

\renewcommand{\eqref}[1]{(\ref{#1})}

\newcommand{\wgt}{\mathrm{wgt}}

\newcommand{\bfs}{{\mathbf s}}
\newcommand{\bfc}{{\mathbf c}}

\newtheorem{prop}{Proposition}[section]

\newtheorem{theo}[prop]{Theorem}
\newtheorem{conj}[prop]{Conjecture}
\newtheorem{exam}{Example}[section]

\begin{document}
\baselineskip 18pt
\title{On Using $(\Z^2, +)$ Homomorphisms to Generate Pairs of Coprime Integers}
\author{Brian A. Benson}

\email{b.a.benson@gatech.edu}

\begin{abstract}
We use the group $(\Z^2,+)$ and two associated homomorphisms, $\tau_0, \tau_1$, to generate all distinct, non-zero pairs of coprime, positive integers which we describe within the context of a binary tree which we denote $T$.  While this idea is related to the Stern-Brocot tree and the map of relatively prime pairs, the parents of an integer pair these trees do not necessarily correspond to the parents of the same integer pair in $T$.  Our main result is a proof that for $x_i \in \{0,1\}$, the sum of the pair $\tau_{x_1}\tau_{x_2} \cdots \tau_{x_n} [1,2]$ is equal to the sum of the pair $\tau_{x_n}\tau_{x_{n-1}} \cdots \tau_{x_1} [1,2]$.  Further, we give a conjecture as to the well-ordering of the sums of these integers.

\end{abstract}

\maketitle

\section{Introduction}
\setcounter{section}{1}

One elementary method of generating all pairs of relatively prime numbers is the map of relatively prime pairs given in \cite{S03}.  In a directed sense, two disjoint roots of the graph are given by $[0,1]$ on the left and $[1,0]$ on the right each with one edge into $[1,1]$ which is located in the center of the graph.  Beginning with the vertex $[1,1]$, each vertex has two children (one to its right and one to its left) whose vertex vector is given by adding the vectors corresponding to its nearest decendent to its left and its nearest decendent to its right.  This process is exactly the same as in the Stern-Brocot tree which, in its entirity, displays all irreducible fractions of $\Q$ sorted by their natural well-ordering from left to right \cite{B07a, GKP90}; specifically, the first entry of the vector in the map of relatively prime pairs gives the numerator while the second entry gives the denomenator of the fraction in the Stern-Brocot tree.  We construct a binary tree of coprime pairs $T$.  While in some ways similar, $T$ differs from the tree given by the map of relatively prime pairs in the sense that two vectors that are neighbors in the map are not necessarily neighbors in $T$.  We must note that $T$ is somewhat related to a binary encoding of the Stern-Brocot tree given in \cite{B07bb}, however, this relationship is not completely explicit and they give no mention of our main result.

Specifically, in this work, we use the $(\Z^2, +)$ and two associated homomorphisms on the element $[1,2] \in \Z^2$ to generate all distinct, non-zero pairs of coprime, positive integers.  This process is given as the explicit construction of $T$.  Specifically, $T$ is generated by two homomorphisms on an element of $(\Z^2, +)$ under ordinary addition.  While this idea is motived by the Cayley graph, $T$ is generated by homomorphisms on a base element, $[1,2]$, of $\Z^2$.  This algebraic construction of $T$ allows us to prove our main result, namely, for $x_i \in \{0,1\}$, the sum of the integers in $\tau_{x_1}\tau_{x_2} \cdots \tau_{x_n} [1,2]$ is the same as the sum of the integers in $\|\tau_{x_n}\tau_{x_{n-1}} \cdots \tau_{x_1} [1,2]\|_1$.  In addition, we give a conjecture corresponding to the natural well-ordering of these sums and their structure on $T$.

\section{Construction of $T$}
\setcounter{equation}{1}

We consider the natural group $(\Z^2, +)$ under ordinary addition; in other words, $[a,b]+[a',b']=[a+a',b+b']$ for $a,a',b,b',c \in \Z$.  Let $[a,b] \in \Z^2$, we define maps $\tau_0 [a,b] = [a,a+b]$ and $\tau_1 [a,b]=[b, a+b]$.

\begin{prop}Let $a,b \in \Z^+$ with $a<b$.  Then, the inverse image of $\tau_0$ and $\tau_1$ is given by $$r([a,b])= \left \{ \begin{array}{lr} \left [a,b -a \right ], & \mbox{ if $b > 2a$}\\ \left [b - a, a \right ], & \mbox{ if $b \leq 2a$.}\\ \end{array} \right. $$  \end{prop}

\Proof If $\tau_0[c,d]=[a,b]$, then $[c,c+d]=[a,b]$ implying that $c=a$ and $d=b-a$.  Further, $\tau_1[c,d]=[a,b]$, then $[d,c+d]=[a,b]$ implying that $c=b-a$ and $d=a$.\eproof

\begin{prop} The maps $\tau_0 , \tau_1$ are homomorphisms on $\Z^2$.\end{prop}

\Proof Let $a,b,a',b' \in \Z^+$.  For $\tau_0$, we have $$\tau_0([a, b])+\tau_0([a', b'])=[a, a+b] + [a', a'+b']= [a+a', a+b+a'+b']=$$ $$\tau_0 ([a+a', b+b'])=\tau_0 ([a,b]+[a', b']).$$ For $\tau_1$, we have $$\tau_1([a,b])+\tau_1([a', b']) = [b,a+b]+[b', a'+b'] = [b+b', a+b+a'+b']=$$ $$\tau_1([a+a', b+b'])=\tau_1([a, b]+[a',b']).$$  \eproof

Note that if we extend $\Z^2$ to a $\Z$-module with action defined as $c[a,b]=[ca,cb]$, we have $c \in \Z$, $\tau_0 [ca,cb]=  [ca,ca+cb] = c [a, a+b]=c \tau_0 [a,b]$ and $\tau_0 [ca,cb]= [cb,ca+cb] = c [b, a+b]=c \tau_1 [a,b]$.  Thus, $\tau_0$ and $\tau_1$ are homomorphisms in the more general setting of a canonical $\Z$-module.

\noindent {\bf Remark:}  The map $r$ is not a homomorphism.  To see this, consider that $[1,4]+[2,3]=[3,7]$, however, $r[1,4]+r[2,3]=[1,3]+[1,2]=[2,5] \neq [3,4]=r[3,7]$.

Now, we construct $T$ by taking the single element $[1,2] \in \Z^2$ and recursively adding an edge from $[1,2]$ to a new vertex corresponding to $\tau_0 [1,2]$ and from $[1,2]$ to a new vertex corresponding to $\tau_1 [1,2]$.  Thus,  $T$ contains the edge between $[a,b]$ and $[a',b']$ if and only if $\tau_0[a,b]=[a',b']$, $\tau_1[a,b] = [a',b']$, or $r [a,b] = [a',b']$.  To see that $T$ is a tree, for the image of $[1,2]$ under $\tau_0$ and $\tau_1$, if $\tau_0[a,b] = \tau_1[c,d]$, then $[a,a+b]=[d,c+d]$; this is a contradiction since, by induction on the repeated images of $\tau_0$ and $\tau_1$ on $[1,2]$, $a<b$ and $c<d$.  Further, since $T$ contains only the elements $[a,b] \in \Z^2$ with $a,b$ both positive and $a<b$, $\tau_i[a,b]=\tau_i[a',b']$ if and only if $a=a'$ and $b=b'$.

\begin{prop} $T$ contains all positive pairs of distinct coprime integers.\end{prop}

\Proof  To see this, we need only take arbitrary $a,b \in \Z^2$ with $\gcd (a,b)=1$ and $a < b$ and show that $[a,b]$ reduces to $[1,2]$ under successive images of $r$.  However, by the definition of $r$, we have $r[a,b] = [a',b']$ where $0<a'<b'$ and $b'<b$ unless $[a,b]=[1,2]$; thus, $[a,b]$ must reduce to $[1,2]$ under successive images of $r$ implying that $[a,b]$ is contained in $T$. \eproof

\section{Sums of Pairs in $T$}
\setcounter{equation}{1}

We now begin to explore additional properties of $T$ by considering the sums of the pairs of integers in $T$.  The main theorem in this work corresponds to the sums of these pairs of coprime integers and the homomorphisms we use to generate them.  Simply stated, this result can be given as the sum of the integers of $\tau_{x_n}\tau_{x_{n-1}}\cdots \tau_{x_1}[1,2]$ is equal to the sum of the integers in $\tau_1 \tau_2 \cdots \tau_{x_n} [1,2]$ for $x_1, \ldots, x_n \in \{0,1\}$.  Purely as a matter of convenience, since the conventional 1-norm of a vector, $\|[a,b] \|_1=|a|+|b|$, is equivalent to the sum of a pair of positive integers, we use it to denote these sums of integers.  As a further matter of convenience, we use $\|\cdot\|$ to denote $\|\cdot\|_1$ since this use does not produce ambiguity herein.  Upon considering some examples of these sums below, we note that the sums of the codes we check are the same as if the code was listed in reverse order.

\begin{exam}$\|T[1011]\|=\|T[1101]\|$\end{exam}
To show this example, we give the computations of $T[1011]$ and $T[1101]$ respectively.  First, $T[1011]$ is given by  $$[1,2]\stackrel{\tau_1}{\mapsto}[2,3]\stackrel{\tau_0}{\mapsto}[2,5]\stackrel{\tau_1}{\mapsto}[5,7]\stackrel{\tau_1}{\mapsto}[7,12]$$ implying that $\|T[1011]=19\|$.  Second, $T[1101]$ corresponds to $$[1,2]\stackrel{\tau_1}{\mapsto}[2,3]\stackrel{\tau_1}{\mapsto}[3,5]\stackrel{\tau_0}{\mapsto}[3,8]\stackrel{\tau_1}{\mapsto}[8,11]$$ implying that $\|T[1101]\|=19$.\esubproof

\begin{exam} $\|T[1010000]\|=\|T[0000101]\|$\end{exam}
Note that $T[1010000]$ corresponds to
\begin{small}$$[1,2]\stackrel{\tau_1}{\mapsto}[2,3]\stackrel{\tau_0}{\mapsto}[2,5]\stackrel{\tau_1}{\mapsto}[5,7] \stackrel{\tau_0}{\mapsto}[5,12]\stackrel{\tau_0}{\mapsto}[5,17]\stackrel{\tau_0}{\mapsto}[5,22]\stackrel{\tau_0}{\mapsto}[5,27]$$\end{small}
while $T[0000101]$ corresponds to
\begin{small}$$[1,2]\stackrel{\tau_0}{\mapsto}[1,3]\stackrel{\tau_0}{\mapsto}[1,4]\stackrel{\tau_0}{\mapsto}[1,5]\stackrel{\tau_0}{\mapsto}[1,6]\stackrel{\tau_1}{\mapsto}[6,7]\stackrel{\tau_0}{\mapsto}[6,13]\stackrel{\tau_1}{\mapsto}[13,19].$$\end{small}  Thus, $\|T[1010000]\|=32=\|T[0000101]\|$.\esubproof

These examples illustrate the idea behind the main theorem.  However, in order to simplify the statement and proof of this result, we introduce some additional notation.

\noindent {\bf Notation:}  Let each $x_i \in \{0,1\}$, we define $\tau_{x_1x_2 \cdots x_n} [a,b] \coloneqq \tau_{x_n}\tau_{x_{n-1}} \cdots \tau_{1} [a,b]$.  Further, in taking $[1,2]$ to be the root of $T$, we can express each integer pair in $T$ by a binary code corresponding to the order of the composition of $\tau_0$ and $\tau_1$ needed to generate the pair; thus, we define $T[x_1x_2\cdots x_n] \coloneqq \tau_{x_1x_2 \cdots x_n} [a,b] = \tau_{x_n}\tau_{x_{n-1}} \cdots \tau_{x_1} [a,b]$.  To further condense our notation, we often denote an arbitrary binary code with the variable $\bfc$; in other words, $\bfc \coloneqq x_1x_2 \cdots x_n$ with each $x_i \in \{0,1\}$.

Now, we define the {\em reflection} of a binary code $\bfc=x_1x_2\cdots x_n$ to be $\refl(\bfc )=x_n x_{n-1} \cdots x_1$.
Using our newly prescribed notation and definition, we can now restate the theorem more concisely.

\begin{theo}
For a binary code $\bfc$, $\|T[\bfc]\| = \|T[\refl(\bfc)]\|$.
\end{theo}

\Proof  The proof is by induction on the length of a binary code.  Clearly, the theorem is true if a given binary code is a palindrome; thereby, we can assume that all codes henceforth are not palindromes.  On a slightly technical note, our proof will require the use of many vectors of $\Z^2$ which do not appear in $T$.

For the pair of coprime integers $[a,b]$ with $a<b$, note that by proposition 2.2, we have $r([a,b])=\sigma[a,b-a]$ where we can let $\sigma$ be the permutation such that $\sigma[a,b-a]=[a,b-a]$ if $2a \leq b$ and $\sigma[a,b-a]=[b-a,a]$ if $2a>b$.

We consider the base case to be binary codes of length 2, in which there are two non-palindrome codes which are are reflections of one another.  These cases are $01$ and $10$ and it is a straightforward computation to show that $\|T[01]\|=7=\|T[10]\|$.  Suppose that the reflection principle is true for codes of integer length $k \leq n-1$ for $n \in \Z_{\geq 2}$.  Now we consider the arbitrary code $x_1x_2 \cdots x_n$ of length $n$.  To proceed, we must consider $\tau_{x_n}(T[x_1x_2 \cdots x_{n-1}])$; from here, we consider two cases, $x_n=0$ and $x_n=1$.

Due to several complications in the case where $x_n=0$, we initially suppose $x_n=1$\footnote{We will be able to use our proof of the simpler case of $x_n=1$ in order to simplify the number of cases that we must consider when $x_n=0$.}; we note that we can represent $\tau_1$ on an integer pair $\bfs=[a,b]$ in $T$ as $\tau_1(\bfs)=\tau_1 ([a,b])=[b,a+b]=[a,b]+[b-a,a]=\bfs+\sigma^{\ast}(r(\bfs))$ up to some permutation $\sigma^{\ast}$\footnote{Note that $\sigma(\bfs)$ arranges the entries of $\bfs$ by well-ordering the integers from least to greatest.  Thus, in general, $\sigma (\bfs ) \neq \sigma^{\ast} (\bfs)$.}.  Therefore, $$T[x_1x_2 \cdots x_n]=\tau_1(T[x_1x_2\cdots x_{n-1}])=T[x_1x_2\cdots x_{n-1}]+\sigma^{\ast}(r(T[x_1x_2\cdots x_{n-1}]))=$$ $$ T[x_1x_2\cdots x_{n-1}]+\sigma^{\ast}(T[x_1x_2\cdots x_{n-2}]).$$  Now, $$T[x_nx_{n-1}\cdots x_2x_1]=\tau_{x_nx_{n-1}\cdots x_1} [1,2] =\tau_{x_{n-1}x_{n-2}\cdots x_1} [2,3]=\tau_{x_{n-1}{x_{n-2}} \cdots x_{1}}([1,1] + [1,2]) =$$ $$\tau_{x_{n-1}x_{n-2}\cdots x_{1}}[1,1] + \tau_{x_{n-1}x_{n-2}\cdots x_{1}}[1,2]=T[x_{n-2}x_{n-3} \cdots x_2 x_1]+T[x_{n-1}x_{n-2} \cdots x_2 x_1].$$  By the induction hypothesis, $\|T[x_{n-2}\cdots x_2x_1]\|=\|T[x_1x_2 \cdots x_{n-2}]\|$ and $\|T[x_{n-1}x_{n-2}\cdots x_1]\|=$

\noindent
$\|T[x_1x_2 \cdots x_{n-1}]\|$.  Now, we have that $$\|T[x_1x_2 \cdots x_n]\|=\|T[x_1x_2 \cdots x_{n-2}]+T[x_1x_2 \cdots x_{n-1}]\|=\|T[x_1x_2 \cdots x_{n-2}]\|+\|T[x_1x_2 \cdots x_{n-1}]\|=$$ $$\|T[x_{n-2}\cdots x_2x_1]\|+\|T[x_{n-1}\cdots x_2x_1]\|=\|T[x_{n-2}\cdots x_2x_1]+T[x_{n-1}\cdots x_2x_1]\|=\|T[x_{n}\cdots x_2x_1]\|$$ proving the case of $x_n=1$.

Now, suppose $x_n=0$\footnote{One of the main difficulties with this case arises from the fact that $\tau_0[0,1]=[0,1]$, so we must alter our argument from the case where $x_n=1$.}, clearly, if $x_1=1$, then the proof of the case where $x_n=1$ will suffice in proving this case as well.  Thus, we can assume that $x_1,x_n=0$.  Without loss of generality, we can assume that $x_{n-k}=1$, $1 \leq k \leq \lfloor n/2 \rfloor$, and for all $j < k$, $x_{n-j} =0$ and $x_j =0$.  To prove the case of $x_n=0$, we will induct on $k$ beginning with the base case of $k=1$.

When $k=1$, we know that $x_1=x_n=0$ and $x_{n-1}=1$.  Then, if $T[x_1x_2 \cdots x_{n-1}]=[a,b]$, then $$T[x_1x_2 \cdots x_n] = T[x_1x_2 \cdots x_{n-1}] + [0,a].$$  Further, $$T[x_nx_{n-1} \cdots x_1]=\tau_{x_{n}x_{n-1}\cdots x_1} [1,2]=\tau_{x_{n-1}x_{n-2} \cdots x_1}[1,2] + [0,1]=$$ $$T[x_{n-1}x_{n-2} \cdots x_1]+\tau_{x_{n-1}x_{n-2} \cdots x_{1}} [0,1].$$

\vspace{2mm}
\noindent
Since $\tau_{x_{n-1}}=1$, we have $$\tau_{x_{n-1} x_{n-2}\cdots x_{1}}([0,1])=\tau_{x_{n-2}x_{n-3} \cdots x_{1} }([1,1])=\tau_{x_{n-3} \cdots x_{1}}([1,2])=T[x_{n-3}x_{n-4} \cdots x_1].$$  Thus, we have $$T[x_nx_{n-1} \cdots x_1]=T[x_{n-1}x_{n-2} \cdots x_1]+T[x_{n-3}x_{n-4} \cdots x_1].$$

Since, by the original induction hypothesis, $\|T[x_1x_2 \cdots x_{n-1}]\|=\|T[x_{n-1}x_{n-2} \cdots x_1]\|$, we must check that $\|T[x_{n-3}x_{n-4} \cdots x_1]\|=a$.  To do this, we again rely on the original induction hypothesis to tell us that $\|T[x_{n-3}x_{n-4} \cdots x_1]\|=\|T[x_1x_2 \cdots x_{n-3}]\|$\footnote{This allows us to compute $\|T[x_1x_2 \cdots x_{n-3}]\|$ in place of the final entry of $\|T[x_{n-3}x_{n-4} \cdots x_1]\|$.}.

Now, since we took $T[x_1x_2 \cdots x_{n-1}]=[a,b]$, we have that $$T[x_1 \cdots x_{n-2}]=r(T[x_1 \cdots x_{n-1}]).$$  Now, $r(T[x_1 \cdots x_{n-1}])= [a,b-a] \text{ or }[b-a,a]$ depending on the ordering of $2a$ and $b$.  However, since we assumed that $x_{n-1}=1$, we consider the equation $$\tau_{x_{n-1}} [a',b'] = \tau_{1} [a',b']=[b',a'+b'] = [a,b]$$ where $a'<b'$.  By our equation, we have that $a=b', b=a'+b'$ which implies the equalities $a'=b-a, b'=a$.  Therefore, $b-a <a<b$ implying that $$T[x_1 \cdots x_{n-2}]=r(T[x_1 \cdots x_{n-1}])= r[a,b]= [b-a,a].$$  Now, $T[x_1 \cdots x_{n-3}]=r(T[x_1 \cdots x_{n-2}])=r[b-a,a]=[|b-2a|,b-a]$ or $[b-a, |b-2a|]$ depending on the ordering of $b-a$ and $|b-2a|$.  Either way, we can conclude that $\|T[x_{n-3} \cdots x_1]\| = \|T[x_1 \cdots x_{n-3}]\| = a$.  Thus, we can conclude when $x_n=0$ and $k=1$, $\|T[x_1 \cdots x_n]\|=\|T[x_n \cdots x_1]\|$ completing the base case of the induction on $k$.

Now, assume that the claim is true for all integers $k$ such that $1<k<m< \lfloor n/2 \rfloor$.  Now, we consider the case of $m+1$.  Then, by this, we know that $x_1, \ldots, x_{m-1}=0$ and $x_{n-m}, \ldots, x_{n}=0$.  If we again take $T[x_1 \cdots x_{n-1}]=[a,b]$, we have that $$T[x_1 \cdots x_n]=T[x_1 \cdots x_{n-1}] + [0,a].$$  Further, $$T[x_n\cdots x_1]=T[x_{n-1} \cdots x_1]+\tau_{x_{n-1}\cdots x_{1}}[0,1].$$  However, since $x_{n-m}, \ldots, x_{n-1}=0$, $\tau_{x_{n-1} \cdots x_{1}}[0,1]=\tau_{x_{n-m+1} \cdots x_{1} }[0,1]=T[x_{n-m-1} \cdots x_1]$.  Similar to the base case, we consider $\|T[x_1 \cdots x_{n-m-1}]\|=$

\noindent
$\|r^m(T[x_1 \cdots x_{n-1}])\|$ which is equal to $\|T[x_{n-m-1}\cdots x_1]\|$ under the induction hypothesis.

Again refering to the pair $[a',b']$ with $a'<b'$, we wish to find the exact permutation $\sigma$ required for $r([a,b])$ where $\tau_0([a',b'])=[a,b]$.  Thus, we consult the equation $$\tau_0 [a',b'] = [a',b'] + [0,a'] = [a,b].$$  Solving for $a'$ and $b'$, we get $a'=a$, $b'=b-a$.  Since $x_{n-m}, \ldots, x_{n}=0$, we know that $\min(T[x_1 \cdots x_{n-m+1}]=\min(r^{m-2}(T[x_1 \cdots x_{n-1}]))=a$.  To find $\|T[x_1 \cdots x_{n-m-1}]\|$, we must find $\|r^2(T[x_1 \cdots x_{n-m+1})\|$; but since $x_{n-m+1}=1$, we have that $$\|r^2(T[x_1 \cdots x_{n-m+1})\| = \min(T[x_1 \cdots x_{n-m+1}] = a$$ by our proof of the base case of the induction on $k$.  Thus, $\|T[x_1 \cdots x_n]\|=\|T[x_1 \cdots x_{n-1}]\|+a = \|T[x_{n-1} \cdots x_{1}]\|+a=\|T[x_n \cdots x_1]\|$ completing the proof of the case $x_n=0$.  Thus, all possible cases of non-palindrome codes have been exhausted, proving the theorem.\eproof

\noindent {\bf Remark.}  The converse of this theorem, for codes of equal length, does not hold true in general.  Consider the codes $10011$ and $01110$; a simple compution yeilds $\|T[10011]\|=25=\|T[01110]\|$ confirming this fact.

Motivated by the theorem, we would like to see how the number of changes or transitions of a code $x_1x_2 \cdots x_k$ affects $\|T[x_1 x_2 \cdots x_k]\|$.  This gives a simple proposition which begins to outline our conjecture.

\begin{prop}
For integer $j \geq 2$, $$\|T[\underbrace{00\cdots 0}_{j}\underbrace{11\cdots 1}_{j}]\|=\|T[\underbrace{11\cdots 1}_{j}\underbrace{00\cdots 0}_{j}]\|<\|T[\underbrace{0101 \cdots 01}_{2j}]\|=T[\underbrace{1010 \cdots 10}_{2j}]\|.$$
\end{prop}

\Proof The equivalent conditions hold trivially by the theorem.  Using this same reasoning, it is only necessary to compare $\|T[\underbrace{0101\cdots01}_{2j}]\|$ and $\|T[\underbrace{11\cdots 1}_{j}\underbrace{00\cdots 0}_{j}]\|$. The base case is $j=2$ where $T[1100]=[3,11,14]$ and $T[0101]=[7,10,17]$; from here we note that $\min(T[1100])<\min(T[0101])$ and $\|T[1100]\|<\|T[0101]\|$.  Now, we apply the step of induction assuming that these two inequality statements are true for all integers $j\leq k$ for arbitrary $k\in \Z^+$.  Since $T[1]=[F_3,F_4]$, where $F_i$ corresponds to the $i$-th integer in the Fibonacci sequence, we note that $T[\underbrace{11\cdots 1}_{j}\underbrace{00\cdots 0}_{j}]=[F_{j+2},F_{j+4}+(j-1)F_{j+2}]$.

Let $T[\underbrace{0101\cdots01}_{2k}]=[a,b]$.  Note that $T[\underbrace{0101\cdots01}_{2(k+1)}]=[a+b,2a+b]$ and $T[\underbrace{11\cdots 1}_{k+1}\underbrace{00\cdots 0}_{k+1}]=[F_{k+3},F_{k+3}+kF_{k+3}]$.  By the induction hypothesis and the fact that the Fibonacci sequence is monotone increasing, we have $a>F_{k+2}>F_{k+1}$.  Thus, $$a+2c=a+2(a+b)>a+2(2a)=5a>2a>2F_{k+2}>F_{k+2}+F_{k+1}=F_{k+3}.$$  By this, we have $$\|T[\underbrace{11\cdots 1}_{k+1}\underbrace{00\cdots 0}_{k+1}]\|=F_{k+5}+(k+1)F_{k+3}=F_{k+4}+F_{k+3}+(k+1)(F_{k+2}+F_{k+1})=$$ $$F_{k+4}+kF_{k+2}+F_{k+3}+F_{k+2}+(k+1)F_{k+1}<a+b+F_{k+3}+F_{k+2}+(k+1)F_{k+1}=$$ $$a+b+F_{k+4}+kF_{k+1}+F_{k+1}<a+b+F_{k+4}+kF_{k+2}+F_{k+2}<3a+2b.$$ This gives the result.\eproof

To give a more generalized notion of this observation, we define several properties of general codes of $T$. Consistent with binary codes in coding theory, we let the weight of a configuration, $\wgt (\bfc)$, represent the number of ones present in the code.  We define the {\em edge cluster number} of the edge $x_i$ in the code $\bfc=x_1x_2\cdots x_i \cdots x_n$ to be $\clus (x_i, \bfc)=|\{x_k:x_k=x_j=x_i \text{ for all } \min(i,k) \leq j \leq \max(i,k)\}|$.  From this, we can define the {\em cluster average} of the code $\bfc$ associated with $T$ to be $$\avg(\bfc)=\sum_{i=1}^{n}{\frac{[\clus(x_i,\bfc)]}{n}}.$$  Further, we define the {\em cluster variance} of the code $\bfc$ to be $$\var(\bfc)=\sum_{i=1}^{n}{\frac{[\clus(x_i,\bfc)]^2}{n}}.$$

The following is an example of the relation between the variance of a code and the sum of its correponding value with respect to $T$.

\begin{exam} Note that $\wgt(\underbrace{0101\cdots 01}_{2j})=\wgt(\underbrace{11\cdots 1}_{j}\underbrace{00\cdots 0}_{j})$ while $\var(\underbrace{0101\cdots01}_{2j})=2j<2j^3=j\cdot j^2+j\cdot j^2=\var(\underbrace{11\cdots 1}_{j}\underbrace{00\cdots 0}_{j})$.  Further, by the proposition, $\|T[\underbrace{11\cdots 1}_{j}\underbrace{00\cdots 0}_{j}]\|< \|F[\underbrace{0101\cdots 01}_{2j}]\|$.\esubproof \end{exam}
\noindent
This example illustrates the beginning of observations which give evidence for the conjecture; its remaining rationale  is given after its statement.

\begin{conj}
For codes $\bfc_1,\bfc_2$ of $T$ of equal length, if $\wgt(\bfc_1)=\wgt(\bfc_2)$ and $\var(\bfc_1)<\var(\bfc_2)$, then $\|T[\bfc_1]\|>\|F[\bfc_2]\|$.
\end{conj}

The overriding rationale behind the conjecture is that the larger the variance of the code, the smaller the average clusters size and, thus, the more transitions there are back and forth between ones and zeros when the weight or the code is constant.  Now, suppose that we wish to maximize $\|T[\bfc]\|$ for a code $\bfc$ of fixed length and weight; then, we suppose that $T[x_1 x_2 \cdots x_k]=[a,b,c]$ which naturally implies that $a,b \in \Z^+$, $a < b$, and $a+b=c$.  Now note that $T[x_1 x_2 \cdots x_k 0]=[a,c,a+c]$ while $T[x_1 x_2 \cdots x_k1]=[b,c,b+c]$ which means that the first entry of $T[x_1 x_2 \cdots x_kx_{k+1}]$ is maximized locally by choosing $x_{k+1}=1$ while the second entry of $T[x_1 x_2 \cdots x_kx_{k+1}]$ is the same irregardless of the value of $x_{k+1}$.  Since the length and weight of $\bfc$ are fixed, when $x_j$ must be zero, having $x_{j-1}=1$ maximizes the sequence locally.  Spreading this local observation over the entire length of the code gives evidence for the conjecture.  Further, although these local observations are relatively straightforward, it appears that constructing a rigorous proof of the conjecture from these observations is somewhat less intuitive.  In the following paragraphs of this section, we will consider a few approaches towards proving the conjecture.

Perhaps the first approach a reader might take is straightforward induction on the code length\footnote{In addition, perhaps even induction on the code weight for each code length as well.}.  However, if $\wgt(x_1 \cdots x_k)=\wgt(x_1' \cdots x_k')$ and $\var [x_1 \cdots x_k] < \var [x_1' \cdots x_k']$, then it is not necessarily true that $\var [x_1 \cdots x_k x_{k+1}] \leq \var [x_1' \cdots x_k'x_{k+1}]$\footnote{Note, however, that if it is the case that $x_k \neq x_{k+1}$, then $\var [x_1 \cdots x_k x_{k+1}]=(n/ (n+1))(\var [x_1 \cdots x_k]+1) < (n/(n+1))(\var [x_1' \cdots x_k']+1) \leq \var [x_1' \cdots x_k'x_{k+1}]$.}.  A counter-example which tells us that this is not true in general is as follows:  $\var[1010111] = 31/7 < 55/7 =\var [1110110]$, however, $\var [10101111] = 17/2 > 65/8 = \var[11101101]$.  

Note that even if the orderings of the variance of $\bfc$ is known with respect to an inductive step, this in itself does not appear to be enough to establish the ratio between the first and second entry of $T[\bfc]$; to determine this ratio from a label without direct computation from the code, it appears that something else must be known about the structure of the code.  However, as the length of the codes under consideration becomes larger, additional structures arise which make this approach non-trivial.

\begin{exam} We consider all codes of length 4 within the context of the conjecture.\end{exam}

The codes of the non-trivial weights are given below.

\begin{center}
\begin{tabular}[t]{ccccc}
 Binary Code $\bfc$ & $\var [\bfc]$ & Generation of $T[\bfc]$ \\
 $1000$ & 7 &  $[1,2] \stackrel{\tau_1}{\mapsto} [2,3] \stackrel{\tau_0}{\mapsto} [2,5] \stackrel{\tau_0}{\mapsto} [2,7] \stackrel{\tau_0}{\mapsto} [2,9]$\\
 $0100$ & 5/2 & $[1,2] \stackrel{\tau_0}{\mapsto} [1,3] \stackrel{\tau_1}{\mapsto} [3,4] \stackrel{\tau_0}{\mapsto} [3,7] \stackrel{\tau_0}{\mapsto} [3,10]$\\
 $0010$ & 5/2 & $[1,2] \stackrel{\tau_0}{\mapsto} [1,3] \stackrel{\tau_0}{\mapsto} [1,4] \stackrel{\tau_1}{\mapsto} [4,5] \stackrel{\tau_0}{\mapsto} [4,9]$\\
 $0001$ & 7 & $[1,2] \stackrel{\tau_0}{\mapsto} [1,3] \stackrel{\tau_0}{\mapsto} [1,4] \stackrel{\tau_0}{\mapsto} [1,5] \stackrel{\tau_1}{\mapsto} [5,6]$\\ 
$1100$ & 4 & $[1,2] \stackrel{\tau_1}{\mapsto} [2,3] \stackrel{\tau_1}{\mapsto} [3,5] \stackrel{\tau_0}{\mapsto} [3,8] \stackrel{\tau_0}{\mapsto} [3,11]$\\
$0110$ & 5/2 & $[1,2] \stackrel{\tau_0}{\mapsto} [1,3] \stackrel{\tau_1}{\mapsto} [3,4] \stackrel{\tau_1}{\mapsto} [4,7] \stackrel{\tau_0}{\mapsto} [4,11]$\\
$0011$ & 4 & $[1,2] \stackrel{\tau_0}{\mapsto} [1,3] \stackrel{\tau_0}{\mapsto} [1,4] \stackrel{\tau_1}{\mapsto} [4,5] \stackrel{\tau_1}{\mapsto} [5,9]$\\
$1001$ & 5/2 & $[1,2]\stackrel{\tau_1}{\mapsto} [2,3] \stackrel{\tau_0}{\mapsto} [2,5] \stackrel{\tau_0}{\mapsto} [2,7] \stackrel{\tau_1}{\mapsto} [7,9]$\\
$1010$ & 1 & $[1,2]\stackrel{\tau_1}{\mapsto} [2,3] \stackrel{\tau_0}{\mapsto} [2,5] \stackrel{\tau_1}{\mapsto} [5,7] \stackrel{\tau_0}{\mapsto} [5,12]$\\
$0101$ & 1 & $[1,2] \stackrel{\tau_0}{\mapsto}[1,3] \stackrel{\tau_1}{\mapsto} [3,4] \stackrel{\tau_0}{\mapsto} [3,7] \stackrel{\tau_1}{\mapsto} [7,10]$\\
$1110$ & 7 & $[1,2] \stackrel{\tau_1}{\mapsto} [2,3] \stackrel{\tau_1}{\mapsto} [3,5] \stackrel{\tau_1}{\mapsto} [5,8] \stackrel{\tau_0}{\mapsto} [5, 13]$\\
$1101$ & 5/2 & $[1,2] \stackrel{\tau_1}{\mapsto} [2,3] \stackrel{\tau_1}{\mapsto} [3,5] \stackrel{\tau_0}{\mapsto}[3,8] \stackrel{\tau_1}{\mapsto} [8,11]$\\
$1011$ & 5/2 & $[1,2] \stackrel{\tau_1}{\mapsto} [2,3] \stackrel{\tau_0}{\mapsto} [2,5] \stackrel{\tau_1}{\mapsto} [5,7] \stackrel{\tau_1}{\mapsto} [7,12]$\\
$0111$ & 7 & $[1,2] \stackrel{\tau_0}{\mapsto} [1,3] \stackrel{\tau_1}{\mapsto} [3,4] \stackrel{\tau_1}{\mapsto} [4,7] \stackrel{\tau_1}{\mapsto} [7,11]$\\
\end{tabular}
\end{center}

\vspace{2mm}
\noindent From this, it becomes clear that the conjecture holds true for all codes of length 4; further, note that the variance alone is not enough to determine the ratio between the first and second entry of $T[\bfc]$.\esubproof

\noindent
{\bf Acknowledgement.} The author would like to thank Ian Fredenberg, Yang Wang, Ernie Croot, Prasad Tetali, and Peter Winkler for their insight and helpful discussions.  A portion of the work herein was conducted while the author attended the REU program at the Georgia Institute of Technology in the summers of 2006 and 2007.

\end{document}